\documentclass{amsart}
\usepackage{amsmath,amssymb,latexsym}
\usepackage{mathpazo}
\usepackage[mathscr]{eucal}
\usepackage{pb-diagram}
\usepackage{mathptmx}
\newcommand{\quant}{\mathfrak{q}}

\newcommand{\calA}{{\mathcal{A}}}

\newcommand{\calB}{{\mathcal{B}}}
\newcommand{\calW}{{\mathcal{W}}}
\newcommand{\calC}{\mathcal{C}}

\newcommand{\scrF}{\mathscr{F}}

\newcommand{\sfG}{\mathsf{G}}
\newcommand{\sfH}{\mathsf{H}}

\newcommand{\C}{\mathbb{C}}
\newcommand{\R}{\mathbb{R}}

\newcommand{\N}{\mathbb{N}}
\newcommand{\Z}{\mathbb{Z}}
\newcommand{\Hom}{\operatorname{Hom}}

\newcommand{\tr}{\operatorname{tr}}

\newcommand{\Tot}{\operatorname{Tot}}

\newcommand{\id}{\operatorname{id}}
\newcommand{\spin}{\mathfrak{sp}}

\newcommand{\Ch}{\operatorname{Ch}}
\newcommand{\Sp}{\operatorname{Sp}}

\newcommand{\cpt}{\text{\tiny \rm cpt}}

\numberwithin{equation}{section}
\theoremstyle{plain}
        \newtheorem{theorem}{Theorem}[section]
        
        \newtheorem{proposition}[theorem]{Proposition}
        \newtheorem{corollary}[theorem]{Corollary}

\theoremstyle{definition}
        
        \newtheorem{remark}[theorem]{Remark}
        \newtheorem{example}[theorem]{Example}
        \newtheorem{question}[theorem]{Question}

\title{On the algebraic index for riemannian \'etale groupoids}
\author{M.J.~Pflaum, H. Posthuma,~\textrm{and} X.~Tang}
\begin{document}
\begin{abstract}
 In this paper we construct an explicit quasi-isomorphism to study the cyclic
 cohomology of a deformation quantization over a riemannian \'etale groupoid.
 Such a quasi-isomorphism allows us to propose a general algebraic index problem
 for riemannian \'etale groupoids. We discuss solutions to that index
 problem when the groupoid is proper or defined by a constant Dirac structure on
 a 3-dim torus.
\end{abstract}
\address{\newline
Markus J. Pflaum, {\tt markus.pflaum@colorado.edu}\newline
         \indent {\rm Department of Mathematics, University of Colorado,
         Boulder, USA}\newline
        Hessel Posthuma,  {\tt H.B.Posthuma@uva.nl}\newline
         \indent {\rm Korteweg-de Vries Institute for Mathematics, University of Amsterdam, The Netherlands} \newline
        Xiang Tang, {\tt xtang@math.wustl.edu}   \newline
         \indent {\rm  Department of Mathematics, Washington University,
         St.~Louis, USA}}
\maketitle

\section*{Introduction}
In \cite{ppt:index} and \cite{ppt:high}, the authors of this paper
studied the algebraic index theory over orbifolds using noncommutative
geometry and deformation quantization. In \cite{ppt:index}, we obtained an
explicit topological formula for the Chern character of an elliptic operator
on a compact riemannian orbifold. With this paper we continue our study of
algebraic index theory over singular spaces modeled by groupoids.
More precisely, omitting the assumption
of properness from \cite{ppt:index} and \cite{ppt:high} we study here
the algebraic index theory of a singular space which can be obtained as the
quotient of an \'etale groupoid equipped with an invariant riemannian metric
on the unit space. Such groupoids are called riemannian \'etale groupoids
and appear naturally in the study of riemannian foliations.

A riemannian foliation \cite{molino:book} is a foliation $(M,\scrF)$ equipped
with a bundle like metric $\eta$. Such a bundle like metric defines a holonomy
invariant metric on the normal bundle of the foliation. Let $X$ be a complete
transversal to the foliation $\scrF$ which means that $X$ is an immersed
submanifold of $M$ which intersects with every leaf of $\scrF$. The holonomy
groupoid $\sfG$ associated to $X$ is an \'etale groupoid and the metric $\eta$
defines a $\sfG$-invariant metric on $X=\sfG_0$. Particularly, if
all the leaves of $\scrF$ are closed in $M$, then the
groupoid $\sfG$ is proper and the quotient space $X/\sfG$ is an
orbifold.

Let $\sfG$ be an \'etale groupoid with an invariant riemannian
metric $\eta$. We consider the convolution algebra
$\calC^\infty\rtimes \sfG$ of compactly supported smooth functions
on $\sfG$. This algebra is well studied in noncommutative geometry.
In particular, its cyclic homology was computed by Brylinsky-Nistor
\cite{bn} and Crainic \cite{crainic}. When there is a
$\sfG$-invariant symplectic form on $\sfG_0$, deformation
quantization of the groupoid algebra $\calC^\infty\rtimes\sfG$ was
constructed by the last author \cite{ta}. As a first step toward the
algebraic index theory, we construct explicit cyclic (co)cycles on
the deformation quantization $\calA^\hbar\rtimes \sfG$ of
$\calC^\infty\rtimes\sfG$. More precisely, we construct a
quasi-isomorphism $\Psi$ from the cyclic chain complex of the
algebra $\calA^\hbar\rtimes \sfG$ to the de Rham complex of
compactly supported simplicial differential forms on the inertia
groupoid $\Lambda\sfG$ of $\sfG$ and a quasi-isomorphism $\Phi$ from
the de Rham complex to the cyclic cochain complex of
$\calA^\hbar\rtimes \sfG$. These quasi-isomorphisms generalize our
constructions in \cite[Sec.~5]{ppt:high}. The new input in this
paper is that instead of working with the inertia orbifold which is
the quotient space defined by the inertia groupoid, we construct
$\Phi$ and $\Psi$ on the simplicial spaces associated to the nerves
of the inertia groupoid $\Lambda \sfG$. This improvement allows to
drop the assumption that the groupoid $\sfG$ under consideration is
proper, and also suggests a further generalization to \'etale
groupoids without invariant riemannian structures.

In algebraic index theory one wants to determine topological formulas for the Chern character
of $K_0$ group elements of a deformation quantization. In the work
\cite{ppt:high}, we obtain a complete answer to this question in the
case of proper \'etale groupoids. To answer this question for
general \'etale groupoids, we need more tools in Lie algebra
cohomology and homological algebras. We plan to discuss this in the
future. Instead, we study an interesting example about a constant
Dirac structure on a 3 dimensional torus, whose idea goes back to
joint work of the third author with Weinstein \cite{tw:dirac}. The
transformation groupoid of the integer $\Z$ acting on a 2
dimensional torus $X$ by irrational rotation appears naturally in
describing the geometry of the characteristic foliation associated
to such a Dirac structure. Furthermore, the Dirac structure defines
a $\Z$ invariant symplectic form on the two dimensional torus $X$.
We discuss the topological formulas of $K_0$ elements of the
deformation quantization $\calA^\hbar_X\rtimes \Z$. Our formulas of
the Chern characters agree with Connes' computations in
\cite{connes:foliation}.

This paper is organized as follows. In Section \ref{Sec:QuantConv},
we discuss the construction of deformation quantization
$\calA^\hbar\rtimes \sfG$ and prove an equivariant normal form
theorem for a $\sfG$-invariant symplectic form on $\sfG_0$; in
Section \ref{Sec:ExpChaMap}, we construct a quasi-isomorphism from
the cyclic chain complex of $\calA^\hbar\rtimes \sfG$ to the
simplicial complex of de Rham complex on the inertia groupoid
$\Lambda \sfG$; and in Section \ref{Sec:IndThm}, we propose a
general index problem and discuss the examples of
proper \'etale groupoids and constant Dirac structures on tori.\\

\noindent{\bf Acknowledgments:}
X. Tang's research is partially supported by NSF grant 0703775.
H. Posthuma acknowledges support by NWO.
M.~Pflaum and H.~Posthuma thank the Department of Mathematics of
Washington University, St. Louis, MO for hosting a research visit.


%
%
\section{Quantized convolution algebras}
\label{Sec:QuantConv}
\subsection{Deformation quantization over \'etale groupoids}
Let $\sfG: \sfG_1 \rightrightarrows \sfG_0$ be an \'etale Hausdorff Lie
groupoid. Moreover, assume that $\sfG_0$ carries a $\sfG$-invariant
riemannian metric $\eta$ and an invariant symplectic form $\omega$.
Both tensor fields can be lifted to invariant forms on $\sfG_1$ via pull-back
by the source or target map. The resulting forms are again denoted by $\omega$
and $\eta$. Note that the Levi-Civita connection $\nabla^\text{\tiny \rm LC}$
corresponding to $\eta$ is  $\sfG$-invariant, and that one can construct
from it a torsion-free $\sfG$-invariant symplectic connection $\nabla$
by putting
\[
  \nabla_Y X : = \nabla^\text{\tiny \rm LC}_Y X + \Delta (X,Y) , \quad
  X,Y \in \Gamma^\infty (T\sfG_1 ) ,
\]
where $\Delta \in \Gamma^\infty (T^*M \otimes T^*M \otimes TM ) $ is
the unique tensor field satisfying $\omega (-, \Delta ) = \Delta'$
with the tensor field
$\Delta' \in \Gamma^\infty (T^*M \otimes T^*M \otimes T^*M ) $ defined by
\[
  \Delta' (X,Y,Z) := \nabla^\text{\tiny \rm LC}_Z \omega (X,Y) +
  \nabla^\text{\tiny \rm LC}_Y \omega (X,Z),
  \quad X,Y,Z \in \Gamma^\infty (T\sfG_1 ) .
\]
One checks easily that $\nabla$ is a torsion-free symplectic
connection on $\sfG_1$ (cf.~\cite[Sec.~2.5]{FedDQIT}).
\begin{example}
  As pointed out in the Introduction, the groupoid associated to a riemannian
  foliation is an \'etale and Hausdorff Lie groupoid which is proper if and
  only if all leaves are closed.
  Other examples of \'etale groupoids are the transformation groupoid
  $\Gamma \ltimes M$ of a discrete group $\Gamma$ acting on a smooth manifold $M$
  or the Haefliger groupoid of local diffeomorphisms on a manifold $M$.
  The cotangent bundle of a (riemannian) \'etale
  groupoid even carries a natural invariant symplectic structure.
\end{example}
Next consider the formal Weyl algebra bundle $\calW \rightarrow \sfG_0$. Its fiber
over a point $x \in \sfG_0$ is given by the formal Weyl algebra over $T_x\sfG_0$
which means by the space of formal power series of the form
\begin{equation}
\label{eq:WeylAlgSec}
  a= \sum_{k\in \N, \alpha \in \N^{2d}} \hbar^k a_{k,\alpha} y^\alpha , \quad
  a_{k,\alpha} \in \C ,
\end{equation}
where $2d = \dim \sfG$ and $y^1,\ldots , y^{2d}$ denotes a symplectic basis
of the tangent space $T_x\sfG_0$. The product on $T_x\sfG_0$ is given by
\begin{equation}
   a \circ b = \left. \left( \exp \left( -\frac{i\hbar}{2} \, \omega^{ij}
   \frac{\partial}{\partial y^i}\frac{\partial}{\partial z^j} \right) \,
   a(y,\hbar) \, b(z,\hbar) \right) \right|_{z=y} .
\end{equation}
Observe that $\calW$ can be identified as the associated bundle
$F\sfG_0 \times_{\spin_d} \mathbb{W} \R^{2d}$, where $F\sfG_0$ denotes the
symplectic frame bundle of $\sfG_0$, $\R^{2d}$ carries its canonical symplectic
structure, and $\mathbb{W} \R^{2d}$ denotes the formal Weyl algebra on
$\R^{2d}$. By assigning degree $1$ to every basis element $y_i$, and degree $2$
to the formal parameter $\hbar$, the formal Weyl algebra bundle becomes a bundle
of filtered algebras. Denote by $\calW_{\geq p}$ the bundle of elements of degree
$\geq p$ which means of all elements of form \eqref{eq:WeylAlgSec} where
$a_{k,\alpha} = 0$, if $2k+ |\alpha| < p$.
Next note that the symplectic connection $\nabla$ lifts to a connection on forms with
values in the Weyl algebra bundle:
\begin{equation}
  \nabla : \Omega^\bullet(\sfG_0, \calW) \rightarrow \Omega^\bullet(\sfG_0, \calW) .
\end{equation}
This connection in general has non-vanishing curvature.

According to Fedosov \cite[Sec.~5.2]{FedDQIT} there exists a flat
connection $D: \Omega^\bullet(M, \calW)$ $\rightarrow$
$\Omega^\bullet(\sfG_0, \calW)$ of the form
\begin{equation}
  D = \nabla + \frac{i}{\hbar}  \left[ A, - \right] ,
\end{equation}
where $A \in \Omega^1 \big( \sfG_0, \calW_{\geq 3}\big)$. By construction, $A$
and $D$ are $\sfG$-invariant since $\omega$ and $\nabla$ are invariant.
Since $D$ satisfies the Leibniz-rule the space of flat sections
\[
  \calA^\hbar \big( \sfG_0 \big) := \big\{ f \in \Gamma^\infty \big( \sfG_0,\calW\big)
  \mid Df =0 \big\}
\]
inherits an associative product $\star$ from $\calW$. Moreover, the symbol map
$\sigma : \calA^\hbar \big( \sfG_0 \big) \rightarrow \calC^\infty \big( \sfG_0 \big)$
which locally is given by
\[
  f = \sum\limits_{k\in\N, \: \alpha \in \N^{2d}}
  f_{k,\alpha}y^\alpha \hbar^k \mapsto f_0 := \sum_{k\in\N} f_{k,0} \hbar^k,
  \quad f_{k,\alpha} \in \calC^\infty \big(\sfG_0\big)
\]
is a linear isomorphism. This implies that
$ \calA^\hbar \big( \sfG_0 \big)$ is a deformation quantization over
$\sfG_0$ which means that the following properties hold true:
\begin{enumerate}
\item
  One has
  \[
   \quant (f_1) \star \quant (f_2) = \sum_{k\in \N} \quant \big( c_k (f_1,f_2) \big) \, \hbar^k
  \]
  for all $f_1,f_2 \in \calC^\infty \big(\sfG_0\big) $. Hereby, $\quant$ is the
  inverse of $\sigma$, and the $c_k$ are appropriate bidifferential operators
  on $\calC^\infty \big(\sfG_0\big) $ such that $c_0$ is the commutative
  product of functions.
\item
  $\quant (1)$ acts as unit element with respect to the product $\star$.
\item
  One has
  \[
    \big[ \quant (f_1), \quant (f_2) \big]_\star := \quant (f_1) \star \quant (f_2) -
    \quant (f_2) \star \quant (f_1) = -i \hbar \quant \big( \{ f_1,f_2\} \big) +
    o \big(\hbar^2\big)
   \]
  for all $f_1,f_2 \in \calC^\infty \big(\sfG_0\big) $.
\end{enumerate}
Note that the bidifferential operators $c_k$ are uniquely determined by
$D$, and that $\star$ is $\sfG$-invariant by construction. Moreover, since $D$
is a differential operator, $\star$ is local and one even obtains a
$\sfG$-sheaf of deformed algebras $\big(\calA^\hbar,\star\big)$.
From this sheaf one can form the crossed product algebra
$ \calA^\hbar \rtimes \sfG $. The underlying linear space is given by
$\Gamma^\infty_\cpt \big( \sfG_1 , s^* \calA^\hbar \big)$, the space of all smooth
functions $F :\sfG_1 \rightarrow \calW$ with compact support such that for every
$g\in \sfG_1$ there is an open neighborhood $U$ and an element
$f\in \calA^\hbar \big( s(U) \big)$ such that $F_{|U}=f\circ s_{|U}$.
The product on $ \calA^\hbar \rtimes \sfG $ is the convolution product given
as follows:
\begin{equation}
  \label{eq:DefConProd}
  F_1 \star_\text{\tiny\rm c} F_2 (g) :=  \sum_{g_1 \, g_2 = g}
  F_1 (g_1)g_2 \star  F_2 (g_2), \quad F_1,F_2 \in
  \Gamma^\infty_\cpt \big( \sfG_1 , s^* \calA^\hbar \big).
\end{equation}
Note that hereby we have used that $\sfG$ acts from the right on the sheaf
$\calA^\hbar$. As pointed out in \cite{ta}, the crossed product algebra
$\calA^\hbar \rtimes \sfG$ forms a deformation quantization of the convolution
algebra $\calC^\infty \rtimes \sfG$ along the noncommutative Poisson structure induced by
$\omega$. Moreover, the invariant algebra $\big( \calA^\hbar \big)^\sfG$ forms a
deformation quantization of the algebra $\calC^\infty (\sfG_0)^\sfG$ of invariant smooth
functions on $\sfG_0$.

Our constructions in this paper rely crucially on the existence of an
invariant riemannian structure. In case this condition fails, a more general
theory has been developed by Bressler, Gorokhovsky, Nest, and Tsygan in
\cite{bgnt}.

For the algebraic index theory over the \'etale groupoid $\sfG$ one has to obtain
a precise understanding of the cyclic cohomology theory of the deformed convolution
algebra $ \calA^\hbar \rtimes \sfG $. It is given as follows with proof provided in
the following section:
\begin{equation}
  \label{eq:CycHomDefConvAlg}
  HC^\bullet \big( \calA^\hbar \rtimes \sfG \big) \cong
  \bigoplus_{k\geq 0} H^{\bullet-2k} \big( B \Lambda \sfG , \C((\hbar)) \big).
\end{equation}
This implies in particularly that the space of traces on
$ \calA^\hbar \rtimes \sfG $ is isomorphic to $H^0 \big(B\Lambda\sfG, \C((\hbar))\big)$.

In Equation \ref{eq:CycHomDefConvAlg}, $B\sfH$ denotes the classifying space of an
\'etale groupoid $\sfH$, and $\Lambda \sfG$ is the inertia groupoid of $\sfG$.
The latter is defined as the
transformation groupoid $B_0 \rtimes \sfG$, where
$B_0 := \{ g \in \sfG_1 \mid s(g) = t(g)\}$
is the space of loops of $\sfG$ on which $\sfG$ acts by conjugation
\[
  \sfG_1 \: {_s\!\times_{\!\sigma_0}} \, B_0 \rightarrow B_0, \quad (h,g) \mapsto hgh^{-1} ,
\]
with $\sigma_0:B_0\to\sfG_0$ the moment map given by $\sigma_0(g)=s(g)$.
Note that, although $B_0$ is in general disconnected,
$\sigma_0$ is an immersion. Observe also that the inertia groupoid $\Lambda\sfG$
is \'etale resp.~proper when $\sfG$ is.

Before we construct the explicit chain map realizing the isomorphism in Equation
\eqref{eq:CycHomDefConvAlg}  let us explain in some more detail how to determine the
cohomology of the classifying space of an \'etale groupoid and its inertia groupoid.
\subsection{Cohomology of the inertia groupoid}
\label{Sec:CohIneGro}
For any \'etale Lie groupoid $\sfG$, let $\sfG_k$ be the space of
composable $k$-tuples of arrows:
\[
  \sfG_k= \big\{(g_1,\ldots,g_k)\in\sfG^k \mid s(g_i)=
  t(g_{i+1}),~i=1,\ldots, k-1\big\}.
\]
These spaces form a simplicial manifold with face operators
$\delta_i:\sfG_k\to\sfG_{k-1}$, $i=1,\ldots, k$ given by
\[
 \delta(g_1,\ldots,g_k)=
 \begin{cases}
 (g_2,\ldots,g_k),& \text{for $i=0$},\\
 (g_1,\ldots,g_ig_{i+1},\ldots, g_k), & \text{for $1\leq i\leq k-1$},\\
 (g_1,\ldots,g_{k-1}),&\text{for $i=k$}.
\end{cases}
\]
With this structure,  one obtains the following differentials on the space of
differential forms on $\sfG_\bullet$. First one has the exterior
derivative
$d:\Omega^p(\sfG_q)\to\Omega^{p+1}(\sfG_q)$. Second, the simplicial
structure defines the differential
$\delta:\Omega^p(\sfG_q)\to \Omega^p(\sfG_{q+1})$ given by
\[
 \delta(\alpha):=\sum_{i=0}^q(-1)^i\delta_i^*(\alpha).
\]
Since $\delta$ is given by pull-backs, the two differentials commute:
$[d,\delta]=0$. The total complex computes the cohomology of the classifying space $B\sfG$:
\[
  H^\bullet\big( \Tot(\Omega^\bullet(\sfG_\bullet),d+\delta\big)
  \cong H^\bullet(B\sfG,\C)
\]

Likewise, one can consider the compactly supported differential forms
$\Omega^\bullet_\cpt (\sfG_\bullet)$. Then one has to push forward
to define the derivative by
\[
\partial(\alpha)=\sum_{i=0}^q(-1)^i(\delta_i)_*(\alpha),
\]
which now {\em lowers} degrees:
$\partial:\Omega^p_\cpt(\sfG_q)\to\Omega^p_\cpt(\sfG_{q-1})$.
One easily verifies that
$\big(\Tot(\Omega_\cpt^\bullet(\sfG_\bullet),d,\partial \big)$
forms a mixed complex.

We will apply this to the classifying space of the inertia groupoid
$\Lambda \sfG = B_0 \rtimes \sfG$.
Rather than using the simplicial space associated to $\Lambda\sfG$, we will use the
isomorphic simplicial space given by the higher Burghelea spaces:
\[
  B_k:=\big\{(g_0,\ldots, g_k)\in \sfG_{k+1} \mid s(g_k)=t(g_0)\big\}.
\]
The face operators are simply given by
\[
\delta_i(g_0,\ldots, g_k)=
\begin{cases}
(g_0,\ldots, g_ig_{i+1},\ldots, g_k),&\text{for $0\leq i\leq k-1$},\\
(g_kg_0,\ldots,g_{k-1}),&\text{for $i=k$}.
\end{cases}
\]
There is in fact a canonical cyclic structure on this space given by
\[
t_k(g_0,\ldots, g_k)=(g_k,g_0,\ldots,g_{k-1}).
\]
We will use this cyclic structure in the following way: suppose that $\mathscr{S}_\bullet$ is a cyclic object in the category $\mathsf{Sh}(\sfG)$ of $\sfG$-sheaves on $\sfG_0$, for example the cyclic sheaf resulting from a sheaf of $\sfG$-algebras. This yields the cyclic
vector space
\[
\Gamma_\cpt\left(B_k,\sigma_k^{-1}\mathscr{S}_k\right),
\]
where  $\sigma_k:B_k\to G_0$ is the map defined by $\sigma_k(g_0,\ldots,g_k)=s(g_0)$.
The cyclic structure is given by combining the cyclic structure on $B_\bullet$ together with the cyclic structure on $\mathscr{S}_\bullet$, where
one twists the $k$-th face operator and the cyclic operator by the automorphism $\theta_{g_0\cdots g_k}$ at $(g_0,\ldots,g_k)\in B_k$.
See \cite{crainic} for explicit formulas.

\subsection{Mixed complexes and S-morphisms}
Finally in this section, let us briefly recall some concepts about mixed complexes
as we will need these later.

Recall that by a \textit{mixed complex} one understands a triple $(C_\bullet,b,B)$
where $C_\bullet$ is a graded object in some abelian category,
$b: C_k \rightarrow C_{k-1}$ is a graded map of degree $-1$ and
$B: C_k \rightarrow C_{k+1}$ a graded map of degree $+1$ such that the relations
$b^2 = B^2 = bB+Bb =0$ are satisfied. A mixed complex gives rise to a first quadrant
double complex $\calB C$
\begin{displaymath}
\begin{diagram}
\node{}\arrow{s,r}{b}
\node{}\arrow{s,r}{b}
\node{}\arrow{s,r}{b}
\node{}\arrow{s,r}{b}
\\
\node{C_3}\arrow{s,r}{b}
\node{C_2}\arrow{s,r}{b}\arrow{w,t}{B}
\node{C_1}\arrow{s,r}{b}\arrow{w,t}{B}
\node{C_0}\arrow{w,t}{B}
\\
\node{C_2}\arrow{s,r}{b}
\node{C_1}\arrow{s,r}{b}\arrow{w,t}{B}
\node{C_0}\arrow{w,t}{B}
\\
\node{C_1}\arrow{s,r}{b}
\node{C_0}\arrow{w,t}{B}
\\
\node{C_0}
\end{diagram}
\end{displaymath}
The \textit{Hochschild homology} $HH_\bullet(C)$ of a mixed complex
$C=(C_\bullet,b,B)$ is defined as the homology
of the $(C_\bullet,b)$-complex. The \textit{cyclic homology} $HC_\bullet(C)$ is
defined as the homology of the total complex associated to the double
complex $\calB C$. In this paper, $C_\bullet $ will always be the Hochschild complex
of a unital algebra, $b$ the Hochschild boundary, and $B$ Connes' coboundary.

Every mixed complex  $(C_\bullet,b,B)$ comes with a natural periodicity morphism
$S: \Tot_\bullet \calB C \rightarrow \Tot_\bullet \calB C$ of degree $-2$ which is given by the
canonical projection  $\Tot_k \calB C \rightarrow \Tot_{k-2}\calB C$.

By an S-morphism between two mixed complexes $(C_\bullet,b_C,B_C)$ and
$(D_\bullet,b_D,B_D)$ one now understands a morphism of complexes
$f: \Tot_\bullet \calB C \rightarrow \Tot_\bullet \calB D$ which commutes with the corresponding
periodicity  maps $S_C$ and $S_D$ (cf.~\cite{loday}). Note that such an
S-morphism $f$ induces a morphism on the corresponding Hochschild complexes.
The important observation which we will silently use throughout this paper then
is that $f$ is a quasi-isomorphism with respect to cyclic homology if and only
if it induces an isomorphism on Hochschild homology
(see \cite[Prop.~2.5.15]{loday}).
\subsection{An equivariant normal form  theorem}
For later purposes in this article we derive an equivariant normal form theorem
out of the following result.
\begin{theorem}[Moser Theorem - Equivariant Version]
\label{thm:moser}
  Let $\sfG$ be a riemannian \'etale groupoid, and $\iota: M \hookrightarrow \sfG_0$ an
  invariant submanifold of the base. Let $\omega_0$ and $\omega_1$ be two $\sfG$-invariant
  symplectic forms on $\sfG_0$ such that ${\omega_0}_x = {\omega_1}_x$ for all $x\in M$.
  Then there exist $\sfG$-invariant open neighborhoods $U_0$ and $U_1$ of $M$ in
  $\sfG_0$ and a $\sfG$-equivariant $\varphi : U_0 \rightarrow U_1$ such that
  $\iota = \varphi \circ \iota$ and $\varphi^* \omega_1 = \omega_0$.
\end{theorem}
\begin{proof}
  We only provide a sketch and extend the proof as presented in \cite[Sec.~7.3]{Sil:LSG}
  to the equivariant setting. To this end first note that the $\sfG$-invariant riemannian
  metric on $\sfG_0$ gives rise to an equivariant exponential map
  $\exp : W \rightarrow \sfG_0$ on a $\sfG_0$-invariant neighborhood $W$ of the zero-section
  of $T\sfG_0$. After restriction of $W$ we can assume that $\exp_{T_x\sfG_0 \cap W}$ is
  injective for all $x\in \sfG_0$.
  Let $N\rightarrow M$ be the normal bundle of  $M$ in $\sfG_0$. Obviously, $\sfG$ acts
  on $N$. The restriction $\psi:= \exp_{|\tilde W} : \tilde W  \rightarrow U_0$ with
  $\tilde W := W \cap N$ and $U_0:= \exp (\tilde W)$ then is a $\sfG$-invariant tubular neighborhood
  of $M \hookrightarrow \sfG_0$. Observe that then $ H: U_0 \times [0,1] \rightarrow U_0 $,
  $(x,t)  \mapsto \psi^{-1} \big( t \psi (x)\big)$ is a $\sfG$-equivariant homotopy
  which by the proof of Poincar\'e's Lemma gives rise to a $\sfG$-invariant $1$-form $\nu$
  on $U_0$ such that $\omega_1 - \omega_0 = d \nu$. After possibly shrinking $U_0$, the
  family of closed $2$-forms $\omega_t= (1-t) \omega_0 + t \omega_1 = \omega_0 + t d\nu$
  is in fact a symplectic family for all $t \in [0,1]$. Now let $V_t$ be the unique vector field on $U_0$ such that
  $\omega_t (V_t , -) = - \nu$.
  Integration of $V_t$ provides (after possibly shrinking $U_0$ further, for a $\sfG$-equivariant
  isotopy $ \varrho: U_0 \times [0,1] \times \sfG_0$ such that $\varrho_t^* \omega_t = \omega_0$ for all
  $t \in [0,1]$. Since ${V_t}_{|M} = 0$, one has ${\varrho_t}_{|M} = \id_M$. By putting
  $\varphi:= \varrho_1$ and $U_1 := \varrho_1 (U_0)$, we obtain the diffeomorphism  with the claimed
  properties.
\end{proof}

  Now assume that $M \subset \sfG_0$ is an invariant submanifold and $\sfG_0$ carries an invariant symplectic
  form $\omega$ such that the restriction $\omega_{|M}$ is a symplectic form on $M$. Observe that by assumptions
  the invariant riemannian and symplectic structures on $\sfG_0$ induce a $\sfG$-invariant almost complex
  structure on $\sfG_0$. This gives rise to a an invariant hermitian structure on the normal
  bundle $\pi: N\rightarrow M$ to $TM$ in $\sfG_0$ which in particular means that $N$ is a complex vector bundle over $M$.
  From the riemannian structure on $\sfG_0$ one can construct a $\sfG$-invariant hermitian connection
  $\nabla^N$ on $N$. Choose local complex bundle coordinates $z$ on $N$, and put
  \begin{equation}
    \label{eq:defnormsymp}
     \omega_N := d \big( z^* \nabla^N z \big) + \pi^* \omega_{|M} .
  \end{equation}
  Then $\omega_N$ is a $\sfG$-invariant symplectic form on a neighborhood of the zero-section of $N$.
  The above equivariant Moser Theorem implies the following.
  \begin{corollary}[Equivariant Normal Form]
  \label{thm:norm}
    There exists a $\sfG$-equivariant diffeomorphism $\varphi: U \rightarrow V$ where $U,V$ are
    invariant sufficiently small open neighborhoods of $M$ in $\sfG$ resp.~the normal bundle $N$
    such that $\varphi^* \omega_N = \omega_{|U}$.
  \end{corollary}


%
%
\section{An explicit chain map}
\label{Sec:ExpChaMap}
In this section we present the computation of the cyclic homology and cohomology of the deformation $\calA^\hbar\rtimes\sfG$ of the convolution algebra of $\sfG$ using an explicit
chain map
\[
  {\rm Tot}\left(\calB C_\bullet\left(\calA^\hbar\rtimes\sfG\right)\right)\stackrel{\Psi}{\longrightarrow}
  {\rm Tot}\left(\Omega^\bullet_c(B_\bullet)\right),
\]
as well as its transpose
\[
{\rm Tot}\left(  \Omega^\bullet(B_\bullet)\right)\stackrel{\Phi}{\longrightarrow}
  {\rm Tot}\left(\calB C_\bullet\left(\calA^\hbar\rtimes\sfG\right)\right).
\]
These maps will implement the isomorphisms in Equation
\ref{eq:CycHomDefConvAlg}.

\subsection{Twisted cyclic densities}
\label{Sec:TwiCycDen}
Observe first that the above defined maps $\sigma_k:B_k\to \sfG_0$
are immersions which embed each connected component as a closed submanifold of
$\sfG_k$ because $\sfG_0$ is equipped with an invariant Riemannian metric. For future
use, we introduce the locally constant function $2\ell:B_k\to\N$ given by the
codimension of the above mentioned embedding.

The pull-back of the Weyl algebra bundle $\sigma_k^*\calW$ is a bundle of unital
associative algebras over $B_k$, that comes equipped with a canonical fiber-wise
family of automorphisms
\[
\theta_{g_0\cdots g_k}\in{\rm Aut}\left((\sigma_k^*\calW )_{(g_0,\ldots,g_k)}\right).
\]
This can be seen as follows: The Weyl algebra bundle $\calW\to\sfG_0$ forms a
$\sfG$-vector bundle, so its fiber $(\sigma_k^*\calW_\sfG)_{(g_0,\ldots,g_k)}=\calW_{t(g_0)}$ will carry an action
of the loop $g_0\cdots g_k$ for every $(g_0,\ldots,g_k)\in B_k$. Since $\sfG$ is \'etale, this loop will
act on the tangent space $T_{t(g_0)}\sfG_0$ and because the symplectic form
$\omega$ is $\sfG$-invariant,
this defines an element $g_0\cdots g_k\in \Sp(T_{t(g_0)}\sfG_0)$. In fact, this defines a section of the bundle of symplectic groups
$\sigma_k^*\Sp(T\sfG_0)\to B_k$ associated to the symplectic vector bundle
$\sigma_k^*T\sfG_0$ over $B_k$. Using the fact that the symplectic group of a symplectic vector space acts by automorphisms on its
associated Weyl algebra, this defines the family of automorphisms $\theta$ above.

Using the invariant metric we obtain a decomposition
\begin{equation}
\label{dec-tang}
T_x\sfG_0=T_xB_k\oplus T_x\sigma_k,
\end{equation}
for $x\in B_k$, where $T\sigma_k$ denotes the normal bundle to $\sigma_k:B_k\to\sfG_0$. Because of the $\sfG$-invariance of the symplectic form $\omega$ on $\sfG_0$, this decomposition is one by
symplectic subspaces. This direct sum decomposition factors the Weyl algebra as the tensor product
$\calW(T_x\sfG_0)=\calW(T_xB_k)\otimes\calW(T_x\sigma_k)$.
Clearly the automorphism $\theta$ is the identity on the first factor.

Let
\[
(\tau_0^{\theta},\tau_2^\theta,\ldots,\tau^\theta_{2d-2\ell})\in {\rm Tot}^{2d-2\ell}\left(\calB C (\calW)\right)
\]
be the twisted cyclic cocycle on the Weyl algebra $(\sigma_k^*\calW)_{(g_0,\ldots,g_k)}=\calW_{s(g_0)}$ as constructed in \cite{ppt:high},
with the twist given by the automorphism $\theta=\theta_{g_0\cdots g_k}$. It has the form
\[
\tau_{2m}^{\theta}=\tau_{2m}\,\#~{\rm tr}_{\theta_g},
\]
where $\#$ is the external product in cyclic cohomology with respect
to the tensor product decomposition
$\calW(T_x\sfG_0)=\calW(T_xB_k)\otimes\calW(T_x\sigma_k)$ induced by
\eqref{dec-tang}. The cocycle $\tau_{2m}$ is the untwisted cocycle
of \cite{ppt:high} on $\calW(T_{t(g_0)}B_k)$ which extends the
Hochschild cocycle of \cite{ffs}, and ${\rm tr}_{\theta_{g_0\cdots
g_k}}$ is the twisted trace of \cite{fe:g-index} on
$\calW(T_{t(g_0)}\sigma_k)$. The definition of this trace uses an
auxiliary invariant almost complex structure  on $\sfG_0$ that we
fix using the invariant Riemannian metric and symplectic form. This
turns $T_x\sigma_k$ into a complex vector space and we have
$\theta_x\in U(T_x\sigma_k)$.

Consider the sheaf $\sigma_k^*\calA^\hbar$ of deformation quantizations over $B_k$.
It has an acyclic resolution
\[
  0\longrightarrow\sigma_k^*\calA^\hbar\longrightarrow\sigma_k^*\calW
  \stackrel{D}{\longrightarrow}\Omega^1_{B_k}\otimes \calW
  \stackrel{D}{\longrightarrow}\ldots
\]
with $D$ a Fedosov connection. The proof of this fact is the same as
that of \cite[Prop.~4.7.]{pptt} given the equivariant normal form in
Cor. \ref{thm:norm}. In a local trivialization of the Weyl algebra
bundle induced by local Darboux coordinates,  we can write
$D=d+\hbar^{-1}[A,-]$, with
$A\in\Omega^1(B_k)\otimes\sigma^*_k\calW$ being invariant. Since two
such trivializations differ by a symplectic transformation, $A$ is
unique up to addition of a $\mathfrak{sp}(T\sfG_0)$-valued one-form.

Denote by $\mathscr{C}_\bullet\left(\sigma_k^{-1}\calA^\hbar\right)=
\sigma_k^{-1}\mathscr{C}_\bullet\left(\calA^\hbar\right)$
the sheaf of Hochschild cochains on the sheaf of algebras
$\sigma_k^{-1}\calA^\hbar$ on $B_k$, equipped with the
twisted Hochschild differential $b_\theta$. The associated cyclic bicomplex
of sheaves has a total complex
\[
\mathscr{BC}_\bullet\left(\sigma_k^{-1}\calA^\hbar\right):=\bigoplus_{2m\leq \bullet}\sigma_p^{-1}\mathscr{C}_{q-2m}\left(\calA^\hbar\right)
.
\]
 It is equipped with the twisted differential $b_\theta+B_\theta$. Let
\begin{equation}
\label{shf-map}
  \psi^{i,k}_{2m}\in\underline{\mathscr{H}om}\left(\mathscr{C}_{2m-i}
  \left(\sigma_k^{-1}\calA^\hbar\right),\Omega^i_{B_k}\right),
\end{equation}
where $\underline{\mathscr{H}om}$ denotes the $\Hom$-sheaf,
be defined by
\begin{equation}
\label{form-shf-map}
\begin{split}
  \psi^{i,k}_{2m}\, & \big( a_0\otimes \ldots \otimes a_{2m-i}  \big):=
  \left(\frac{1}{\hbar}\right)^i\tau_{2m}^\theta
  \big( (a_0 \otimes \ldots \otimes a_{2m-i}) \times
  (\sigma_k^* A )_i \big),
\end{split}
\end{equation}
where $\times$ means the shuffle product on Hochschild chains as in
\cite[\S 3.2]{ppt:high}.
\begin{proposition}
\label{prop-shf-map}
The $\psi^{i,k}_{2m}$ satisfy
\[
\begin{split}
(-1)^id\psi^{i,k}_{2m}&=\psi^{i,k}_{2m}\circ b_\theta+\psi^{i,k}_{2m+2}\circ B_\theta,\\
\delta_i^*\psi^{i,k-1}_{2m}&=\psi^{i,k}_{2m},\quad\mbox{\rm for all}~ 0\leq i\leq k.
\end{split}
\]
\end{proposition}
\begin{proof}
First remark that the map $\sigma_k:B_k\to\sfG_0$ factors as
\[
\sigma_k=\sigma_0\circ\lambda_k,
\]
where $\lambda_k:B_k\to B_0$ is given by $\lambda_k(g_0,\ldots,g_k)=g_0\cdots g_k$.
While $\sigma_0$ is an immersion, $\lambda_k$ is \'etale as it is nothing but the $k$-fold
composition of the face maps $\delta_0$.
For $k=0$, the first equation reads
\[
(-1)^id\circ \psi^{i,0}_{2m}=\psi^{i,0}_{2m}\circ b_\theta+\psi^{i,0}_{2m+2}B_\theta,
\]
and was proved in \cite[Prop. 5.5]{ppt:high}. Using the factorization of $\sigma_k$ above
and the fact that $\lambda_k$ is \'etale, the general case follows immediately.

As for the second, this should be interpreted as follows: obviously, the sheaf of differential forms is
simplicial, i.e., there are canonical isomorphisms $\delta_i^*\Omega^p_{B_{k-1}}\cong\Omega^p_{B_k}$
for all $0\leq i\leq k$. The same holds true for $\sigma_p^{-1}\calA^\hbar$,
where one has to be careful with the last one, i.e., $\delta_k$ as it involves
the action of $\sfG$ on $\calA^\hbar$. With this, the sheaf
$\mathscr{C}_{2m-i}\left(\sigma_\bullet^{-1}\calA^\hbar\right)$,
and therefore the $\Hom$-sheaf in \eqref{shf-map} is simplicial. Property $ii)$ now says that
the $\psi^{i,\bullet}_{2m}$ are compatible with all these isomorphisms of sheaves
obtained by pull-back along the face maps.
But this is easy to verify, given the factorization of $\sigma_k$ and the fact that $\lambda_k$ is compatible with
the structure maps.
In particular, for $\delta_k$ we use that the connection $A$ is $\sfG$-invariant.
\end{proof}
On the total sheaf of cyclic chains $\mathscr{BC}\left(\sigma_k^{-1}\calA^\hbar\right)$, this proposition immediately implies:
\begin{corollary}
\label{Cor:SMor}
The morphism
\[
\psi^p:=\sum_{m=0}^{2d-2\ell} \psi^{2d-2\ell-\bullet,p}_{2d-2\ell-2m}:
\left(\sigma_p^{-1}\mathscr{BC}_\bullet\left(\calA^\hbar\right),b_\theta,
B_\theta\right)\to \left(\Omega^{2d-2\ell-\bullet}_{B_p},d,0\right)
\]
is an S-morphism of mixed complexes of sheaves.
\end{corollary}
Recall that $\ell$ is the locally constant function on $B_p$ measuring the codimension in $\sfG$. The formula in Corollary \ref{Cor:SMor} is therefore to be understood as a separate one on each of the connected components of $B_p$. With this, $2d-2\ell$ is by definition the dimension of the connected component of $B_p$ under consideration.
Remark that with the grading as on the right hand side, the de Rham differential {\em lowers} degree,
so that the right hand side indeed forms a mixed complex. The degree zero part of the $S$-morphism of
complexes of sheaves is given by
\[
(\psi^p)_0=\psi^{2d-2\ell-\bullet,p}_{2d-2\ell}:\left(\mathscr{C}_\bullet\left(\sigma_p^{-1}\calA^\hbar\right),b_\theta\right)
\to
\left(\Omega^{2d-2\ell-\bullet}_{B_p},d\right),
\]
explicitly given by formula \eqref{form-shf-map}. Since $B_p$ is symplectic, we can consider Darboux coordinates $(x_1,\ldots,x_{2d-2\ell})$
in a local chart. In such coordinates, the Fedosov connection is gauge equivalent to
\[
\sigma_p^*A=\sum_{i,j=0}^{2d-2\ell}\sigma_p^*\omega_{ij}y_idx_j,
\]
with $y_i$ the element in the Weyl algebra corresponding to $x_i$. Using this in \eqref{form-shf-map} together with the explicit formula
for $\tau_{2d-2\ell}$ in \cite[\S 2]{ppt:high}, one finds
\[
\psi^{2d-2\ell-i,p}_{2d-2\ell}(a_0\otimes\ldots\otimes a_i)=*(a_0da_1\wedge\ldots\wedge da_i)\quad\mbox{mod}~\hbar,
\]
where $*$ denotes the symplectic Hodge star operator $*:\Omega_{B_p}^\bullet\to\Omega_{B_p}^{2d-2\ell-\bullet}$ on $B_p$
as introduced in \cite{bry}. It is known that the quasi-classical limit of the Hochschild
differential is given by
\[
\begin{split}
b(a_0\otimes\ldots\otimes a_k)=&\sum_{i=0}^{k-1}(-1)^i a_0\otimes\ldots\otimes\{a_i,a_{i+1}\}\otimes
\ldots\otimes a_k\\
&+(-1)^k\{\theta(a_k),a_0\}\otimes a_1\otimes \ldots\otimes a_{k-1}\quad \mbox{mod}~\hbar,
\end{split}
\]
which indeed is a complex by the fact that $\theta$ acts by Poisson
automorphism, because the symplectic form on $\sfG_0$ is
$\sfG$-invariant. This is the twisted Poisson homology complex on
$B_k$ of the sheaf $\sigma^{-1}_p\calC^\infty_{\sfG_0}$ as in
\cite[\S 4]{nppt}. The well-known morphism
$\sigma_k^{-1}\mathscr{C}_\bullet\left(\calC^\infty_{\sfG_0}\right)\to\sigma^{-1}_k\Omega^\bullet_{\sfG_0}$
given by
\begin{equation}
\label{hkr}
a_0\otimes\ldots\otimes a_k\mapsto a_0da_1\wedge\ldots\wedge da_k
\end{equation}
defines a quasi-isomorphism to the more familiar Poisson complex on
differential forms, however with the differential similarly twisted
by $\theta$ as in \cite{bn,crainic}. Applying the equivariant normal
form theorem Cor. \ref{thm:norm}, one can prove in the same way as
Fedosov \cite[Sec. 5]{fe:g-index} that the sheaf
$\sigma^{-1}_p\calA^\hbar_{\sfG_0}$ can be identified as the sheaf
of deformation quantization of Weyl algebra valued functions on
$B_p$. With such an observation, one argue in the same way as in
\cite[Prop 5.4.]{nppt} that there exists an inclusion
$\phi_0:\calC^\infty_{B_p}\hookrightarrow\sigma_p^{-1}\calC^\infty_{\sfG_0}$
of sheaves of Poisson algebras which induces a quasi-isomorphism on
the associated Poisson homology complexes and extends to a sheaf
morphism
$\phi:\calA^\hbar_{B_p}\hookrightarrow\sigma_p^{-1}\calA^\hbar_{\sfG_0}$
where $\calA^\hbar_{B_p}$ is an appropriate deformation quantization
of the symplectic manifold $(B_p,\sigma_p^*\omega)$. Because in the
semiclassical limit $\phi_0$ induces an isomorphism on Poisson
homology, $\phi$ is a quasi-isomorphism on the level of Hochschild
chains. We now precompose $\psi^{2d-2\ell-i,p}_{2d-2\ell}$ with
$\phi$.  Since the $*$-operator maps the Poisson homology
differential to the de Rham differential under $*$, one easily
observes that in the quasi-classical limit of
$\phi\circ\psi^{2d-2\ell-i,p}_{2d-2\ell}$ the twisting in the normal
directions drops out, and reduces to \eqref{hkr} composed with $*$
on the symplectic manifold $B_p$. But in this -"untwisted"- case,
this is known to be a quasi-isomorphism, cf.\ \cite{bry}, by the
usual spectral sequence argument. Finally it follows that the
morphism $\psi^{2d-2\ell-i,p}_{2d-2\ell}$ on the Hochschild complex
is a quasi-isomorphism. Since this is just the degree zero part of a
full $S$-morphism $\psi^p$, the $SBI$ sequence implies:
\begin{proposition}
The $S$-morphism $\psi^p$ is a quasi-isomorphism.
\end{proposition}
\subsection{The chain map}
\label{Sec:ChaMap}
With the simplicial morphism $\psi^\bullet$ of mixed complexes of sheaves on
$B_\bullet$, we can now prove our main theorem on the cyclic theory of the
deformed convolution algebra $\calA^\hbar\rtimes\sfG$ of the \'etale groupoid
$\sfG$:
\begin{theorem}
\label{chain-mor}
There exists a natural $S$-morphism
\[
\Psi:\left(\Tot_\bullet \left(\calB C_\bullet\left(\calA^\hbar\rtimes\sfG\right)\right),b,B\right)
\rightarrow
\left(\Tot_\bullet\left(\Omega^{2d-2\ell-\bullet}_c\left(B_\bullet,\C((\hbar))\right)\right),d+
\partial,0\right)
\]
which induces isomorphisms
\[
\begin{split}
HH_\bullet\left(\calA^\hbar\rtimes\sfG\right)&\cong H_\cpt^{2d-2\ell-\bullet}\left(\Lambda\sfG,\C((\hbar))\right),\\
HC_\bullet\left(\calA^\hbar\rtimes\sfG\right)&\cong\bigoplus_{k\geq 0} H_\cpt^{2d-2\ell+2k-\bullet}\left(\Lambda\sfG,\C((\hbar))\right).
\end{split}
\]
\end{theorem}
\begin{proof}
The morphism $\Psi$ is the composition of the following three maps:
\begin{itemize}
\item[$i)$] {\em Restriction to loops}: Since we have $C_k\left(\calA^\hbar\rtimes \sfG\right)\cong\calA^\hbar_{\cpt}\left(\sfG^{\times(k+1)}\right)$, restriction to $B_k\subset\sfG^{\times(k+1)}$ induces a map
\[
  C_k\big( \calA^\hbar\rtimes \sfG \big)\to
  \Gamma_\cpt\left(B_k,\mathscr{C}_k\left(\sigma_k^{-1}\calA^\hbar\right)\right).
\]
As in \cite{bn,crainic} this map constitutes a morphism of cyclic objects where
the right hand side carries the cyclic structure defined by the one on
$B_\bullet$ combined with twisted cyclic structure on
$\sigma_\bullet^{-1}\mathscr{C}_\bullet\left(\calA^\hbar\right)$.
More explicitly, this means that the Hochschild differential
\[
  b:\Gamma_\cpt\left(B_k,\mathscr{C}_k\left(\sigma_k^{-1}\calA^\hbar\right)\right)
  \to \Gamma_\cpt\left(B_{k-1},
  \mathscr{C}_{k-1}\left(\sigma_{k-1}^{-1}\calA^\hbar\right)\right)
\]
is given by applying $\delta_*$ on the domain of sections combined with the
stalk-wise differential $b_\theta$ on the image.
\item[$ii)$] The cyclic {\em Alexander--Whitney} map: consider the bigraded
vector space
\[
  \left(\Gamma_\cpt\big(B_p,\mathscr{C}_q\big(\sigma_p^{-1}\calA^\hbar\big)\big)
  \right)_{p,q\geq 0}.
\]
It is a cylindrical vector space, c.f.~\cite{gj}, with the simplicial and cyclic
operators in the $p$-direction are given by those of the cyclic manifold
$B_\bullet$, whereas in the $q$-direction they are given by the $\theta$-twisted
cyclic structure associated to the stalks of the sheaf of algebra
$\sigma_p^{-1}\calA_\sfG$. As a bisimplicial space, the Eilenberg--Zilber theorem
states that the Alexander--Whitney map
\[
\Gamma_\cpt\left(B_k,\mathscr{C}_k\left(\sigma_k^{-1}\calA^\hbar\right)\right)
\to
\bigoplus_{p+q=k}\Gamma_\cpt\left(B_p,\mathscr{C}_q\left(\sigma_p^{-1}\calA^\hbar\right)\right)
\]
is a quasi-isomorphism from the diagonal to the total simplicial space. Unfortunately, it is not
a morphism of cyclic modules, as one can easily check from the explicit formulas \cite[\S 8.5]{weibel}.
However, the Alexander--Whitney map is the degree zero part of a $S$-morphism which therefore
is a quasi-isomorphism as well. In fact, applying the homological perturbation Lemma, one can get explicit
formulas for this map \cite{kr}.
\item[$iii)$] {\em The map} $\psi$ in \eqref{shf-map}. By Proposition \ref{prop-shf-map} $i)$, the morphism of sheaves
\[
  \sum_{2m\leq q} \psi^{q,p}_{2q-2m}:\bigoplus_{2m\leq q}\sigma_p^{-1}\mathscr{C}_{q-2m}
  \left(\calA^\hbar\right)\to \Omega^q_{B_p}
\]
 is an S-morphism of mixed complexes of sheaves, intertwining the differential $b_\theta+B_\theta$ with $d$. By \ref{prop-shf-map} $ii)$, it is
simplicial over $B_\bullet$.
\end{itemize}
Combining these three maps, one finds the morphism $\Psi$. Since each of its three constituent maps are quasi-isomorphisms as explained
above, the theorem follows.
\end{proof}
As any $S$-morphism, the chain morphism above can be decomposed according to degree. Since in $ii)$ above, the degree zero part is exactly the Alexander--Whitney map, whose explicit for is given e.g. we find:
\begin{corollary}
On the level of Hochschild homology, the isomorphism above is induced by chain map
\[
\Psi_0:C_\bullet\left(\calA^\hbar\rtimes\sfG\right)\to\bigoplus_{p+q=\bullet}\Omega^{2d-2\ell-p}\left(B_q,\C((\hbar))\right)
\]
given by
\[
\Psi_0(a_0\otimes\ldots\otimes a_k)=\sum_{p+q=k}\psi^{2d-2\ell-p,q}_{2d-2\ell}\left(\left.a_0\star\cdots \star a_q\otimes a_{q+1}\otimes\ldots\otimes a_k\right|_{B_q}\right).
\]
\end{corollary}
In the formula above, the products are given by the point-wise, noncommutative multiplication $\star$
in $\calA^\hbar$.
\subsection{The dual map}
\label{Sec:DuaMap}
The dual map
\[
\Phi:\left({\rm Tot}^\bullet\left(\Omega^\bullet(B_\bullet\right),d+\delta,0\right)\to\left({\rm Tot}^\bullet\left(\calB C_\bullet\left(\calA^\hbar\rtimes\sfG\right)\right),b,B\right)
\]
follows by integration: for $\alpha\in\Omega^p(B_q)$ we have
\[
\Phi(\alpha)(a_0\otimes\ldots\otimes a_{p+q}):=\int_{B_q}\alpha\wedge\Psi^{p,q}(a_0\otimes\ldots\otimes a_{p+q}),
\]
where $\Psi^{p,q}:C_{p+q}\left(\calA^\hbar\rtimes\sfG\right)\to\Omega^{2d-2\ell-p}_c\left(B_q,\C((\hbar))\right)$ are the components
of the chain morphism $\Psi$ of Theorem \ref{chain-mor}.

It is easily verified that $\Phi$ is a chain morphism: let $\alpha=(\alpha_{p,q})_{p+q=k}$
be a simplicial differential form on $B_\bullet$ of degree $k$, with
$\alpha_{p,q}\in\Omega^p(B_q)$. Then we have
\[
\begin{split}
\Phi((d+\delta)\alpha)(a)&=\sum_{p+q=k}\int_{B_q}d\alpha_{p,q}\wedge\Psi^{p-1,q}(a)
+\int_{B_{q+1}}\delta\alpha_{p,q}\wedge\Psi^{p,q}(a)\\
&=\sum_{p+q=k}-\int_{B_q}\alpha_{p,q}\wedge d\Psi^{p-1,q}(a)+\int_{B_{q}}\alpha_{p,q}\wedge\partial\Psi^{p,q}(a)\\
&=\sum_{p+q=k}\int_{B_q}\alpha_{p,q}\wedge\Psi^{p,q}((b+B)(a))\\
&=((b+B)\Phi(\alpha))(a),
\end{split}
\]
where $a\in\Tot_k\left(\calB C_\bullet\left(\calA^\hbar\rtimes\sfG\right)\right)$. It
follows from Theorem \ref{chain-mor} that $\Phi$ is an $S$-morphism, i.e., compatible with
the $S$-maps, and induces isomorphisms
\[
\begin{split}
  HH^\bullet\left(\calA^\hbar\rtimes\sfG\right)&\cong H^{\bullet}\left(\Lambda\sfG,\C((\hbar))\right),\\
  HC^\bullet\left(\calA^\hbar\rtimes\sfG\right)&\cong\bigoplus_{k\geq 0} H^{\bullet-2k}\left(\Lambda\sfG,\C((\hbar))\right).
\end{split}
\]



\section{An index problem and examples}
\label{Sec:IndThm}
In this section, we discuss the applications of the
quasi-isomorphism $\Phi$ constructed in the previous section to
algebraic index problem. In Section \ref{Sec:AlgIndPair}, we will briefly recall our
general set up of algebraic index theory; in \ref{Sec:IndOrb}, we will
explain an index theorem in the case that the riemannian \'etale groupoid is
proper; in Section \ref{Sec:IndConstDirac}, we will discuss an index theorem
for quantization of a constant Dirac structure on a 3-dim torus.
\subsection{Algebraic index pairing}
\label{Sec:AlgIndPair}
We recall the general theory of the noncommutative Chern character. Let $A$ be a unital algebra, and $M_n(A)$ be the algebra of $n\times n$ matrices with coefficient in $A$. Idempotent elements in $M_n(A)$ for $n\in \mathbb{N}$ modulo equivalent relations comprise the $K_0$-group of $A$. Let $I$ be the identity matrix in $M_n(A)$. The $k$-th Chern character of a $K_0$ element of $A$ represented by the idempotent $e\in M_n(A)$ is defined as follows,
\begin{equation}\label{eq:chern}
\begin{split}
Ch_k(e)&:=(c_k, \cdots, c_0)\in \calB C_{2k}(A)\\
c_i&:=(-1)^i\frac{(2i)!}{i!}\sum_{s_0, \cdots, s_i}(e-\frac{1}{2}I)_{s_0s_1}\otimes e_{s_1s_2}\otimes \cdots e_{s_{2i}s_0}\in A\otimes \overline{A}^{\otimes 2i},\qquad i\geq 1\\
c_0&:=\sum_{s_0}e_{s_0s_0}\in A.
\end{split}
\end{equation}
It can be easily checked that $Ch_k(e)$ is a cyclic cycle in the
normalized $(b,B)$ bicomplex $\calB C(A)$, and furthermore that equivalent
idempotents in $K_0(A)$ give rise to a same homology class
in $HC_{2k}(A)$. Therefore, $Ch_k$ defines a map from $K_0(A)$ to
the $2k$-th cyclic homology group $HC_{2k}(A)$.

By definition, there is a natural pairing between the cyclic
homology and cohomology groups of an algebra $A$. Given a degree
$2k$ cyclic cocycle $\phi$ of $A$, the index pairing between a
$K_0$-group element $e$ and the cyclic cocycle $\phi$ is defined to
be
\begin{equation}\label{eq:ind-pairing}
\langle e, \phi \rangle:=\langle Ch_k(e), \phi \rangle.
\end{equation}

We remark that if $A$ is not unital, then we need to adjoin a unit to $A$ in order to use the formula (\ref{eq:chern}) to define the Chern character $Ch_k$. For an element $e$ in $K_0(A)$, $Ch_k(e)$ is again a well defined class in $H_{2k}(A)$, and the index pairing defined by Eq. \eqref{eq:ind-pairing} naturally extends.

Applying this index pairing (\ref{eq:ind-pairing}) to the algebra
$\calA^\hbar\rtimes \sfG$, we come up with the following natural
question.
\begin{question}Let $\sfG$ be a Hausdorff \'etale groupoid equipped with an
invariant riemannian metric on $\sfG_0$. Assume that there is a
$\sfG$ invariant symplectic form $\omega$ on $\sfG_0$. We consider the
quantization $\calA^\hbar\rtimes \sfG$ as constructed in Section \ref{Sec:QuantConv}, and the $S$-morphism $\Psi:\big({\rm Tot}_\bullet (\calB C(\calA^\hbar\rtimes\sfG)),b,B\big)
\rightarrow
\big({\rm Tot}_\bullet\big(\Omega^{2d-2\ell-\bullet}_c\big(B_\bullet,\C((\hbar)) \big)\big),d+\partial,0\big)$ defined in Theorem \ref{chain-mor}.
Let $e$ be an element in $K_0(\calA^\hbar\rtimes \sfG)$. What is the
cohomological formula for the class $\Psi(Ch_k(e))$?
\end{question}

We need to develop some Lie algebra cohomology tools like \cite[Sec. 5]{ppt:index} in order to have a full answer to the above question. We plan to address this issue in future publications. In the following two subsections, we will discuss the answer to
this question in two special cases.
\subsection{Index theorem for orbifolds}
\label{Sec:IndOrb} In this subsection, we assume the groupoid $\sfG$
to be proper \'etale. This implies that the quotient space
$X=\sfG_0/\sfG$ is an orbifold. In this case, the cyclic homology
and cohomology of the algebra $\calA^\hbar\rtimes \sfG$ were
computed in full generality in \cite[Sec. 5]{nppt}, and the index pairing
(\ref{eq:ind-pairing}) on the algebra $\calA^\hbar\rtimes \sfG$  was
computed in \cite[Sec. 5]{ppt:high}. We explain  the results in
\cite{ppt:high} briefly.

We start with the observation that the $K_0$-group of an algebra is
invariant under deformations. Therefore, the group $K_0
\big(\calA^\hbar \rtimes \sfG\big)$ is isomorphic to $K_0
\big(\calC^\infty \rtimes \sfG\big)$. According to
\cite{philips:k-theory}, the $K_0$-group of the $C^*$-algebra
completion of $\calC^\infty \rtimes \sfG$ is computed by the
topological $K_0$-group of $\sfG$-equivariant vector bundles on
$\sfG_0$ since $\sfG$ is a proper \'etale groupoid. This group is also called the orbifold $K$-theory group of the underlying orbifold $X$. An element of
topological $K_0$-group of $\sfG$-equivariant vector bundles
consists of a pair of two $\sfG$-equivariant vector bundles $E$ and
$F$ on $\sfG_0$, which are isomorphic outside a compact subset.
Given such a pair of equivariant bundles $(E, F)$, we define an
element $P(E)-P(F)$ in the group $K_0 \big( \calC^{\infty}\rtimes \sfG \big)$,
where $P(E)$ and $P(F)$ are idempotents in
$\calC^{\infty,+}\rtimes \sfG$ and $P(E)-P(F)$ is compactly
supported. (We have used $\calC^{\infty, +}\rtimes \sfG$ for the
groupoid convolution algebra with a unit adjoined.) We refer to
\cite[Sec. 1]{ppt:index} for the details of this construction.
Furthermore, using the standard trick in \cite[Sec.~6.1]{FedDQIT}, we
obtain an element $\widehat{P(E)}-\widehat{P(F)}$ in
$K_0(\calA^\hbar \rtimes \sfG)$ with the similar properties.

To compute the index pairing (\ref{eq:ind-pairing}), we look at the cyclic
cohomology of $\calA^\hbar\rtimes \sfG$. As is explained in Section
\ref{Sec:DuaMap}, the map
\[
  \Phi:\Tot^\bullet\left(\Omega^\bullet(B_\bullet\right))\to
  \Tot^\bullet\left(\calB C\left(\calA^\hbar\rtimes \mathsf{G}\right)\right)
\]
is a quasi-isomorphism. Therefore, the cyclic cohomology of
$\calA^\hbar\rtimes \sfG$ is isomorphic to the cohomology of ${\rm
Tot}^\bullet(\Omega^\bullet(B_\bullet))$. As $\sfG$ is proper, the
groupoid sheaf cohomology on $\sfG$ is zero except at degree 0. And
the degree 0 cohomology of $\sfG$ is equal to the space of invariant
sections. This shows the cohomology of ${\rm
Tot}^\bullet(\Omega^\bullet(B_\bullet))$ is equal to the cohomology
of $\Tot^\bullet(\Omega^\bullet(B_0)^\sfG)= \Tot^\bullet(\Omega^\bullet(B_0/\sfG))$, which is equal to the
cohomology of the quotient $B_0/\sfG$. The space $B_0/\sfG$ is usually called the inertia orbifold $\tilde{X}$ associated to the orbifold $X=\sfG_0/\sfG$.

In order to state our theorem, we introduce characteristic classes
for a $\sfG$-equivariant vector bundle $V$ on $B_0$. Let $R_V$ be
the curvature form of a connection on the bundle $V$ on $B_0$.
Define $Ch_\theta(V)\in H^\bullet(\widetilde{X})$ by
\[
Ch_\theta(V):=\tr\left(\theta \exp(\frac{R_V}{2\pi \sqrt{-1}})\right).
\]
We consider the normal bundle $N$ of $B_0$ embedded as submanifold
of $\sfG$. It is easy to check that $N$ is $\sfG$ equivariant.
Define $Ch_\theta(\lambda_{-1}N)\in
H^\bullet(\widetilde{X})$ by
\[
Ch_\theta(\lambda_{-1}N):=\sum (-1)^\bullet Ch_\theta(\wedge^\bullet N).
\]

We proved the following theorem in \cite{ppt:high}.
\begin{theorem}\label{thm:orb-alg-index} Let $\sfG$ be a proper \'etale Lie groupoid. Consider $E$ and $F$ be a pair of $\sfG$ equivariant vector bundles on $\sfG_0$ defining a $K_0$-group element $\widehat{P(E)}-\widehat{P(E)}$ in $K_0(\calA^\hbar\rtimes \sfG)$, and  $\alpha=(\alpha_{2k}, \cdots, \alpha_0)\in \Tot^{2k}
\Omega^\bullet(\widetilde{X}) ((\hbar))$ a sequence of closed forms. Then we have
\begin{eqnarray*}
  &&\left\langle \Phi(\alpha), Ch_k(\widehat{P(E)}-\widehat{P(F}))\right\rangle\\
  &&\qquad=
  \sum_{j=0}^{k}\int_{\widetilde{X}} \frac{1}{(2\pi \sqrt{-1})^j m}\frac{\alpha_{2j} \wedge
  \hat{A}(B_0/\sfG) \, \Ch_\theta(\iota^*E-\iota^*F) \,
  \exp(-\frac{\iota^*\Omega}{2\pi \sqrt{-1}\hbar})}{\Ch_\theta(\lambda_{-1}N)},
\end{eqnarray*}
where $\iota^*E$ and $\iota^*F$ are pullbacks of $E$ and $F$ to $B_0$ along the source map (same as the target map), and $m$ is a local constant function defined by the order of the isotopy group of the principal stratum of a sector in $\widetilde{X}$, and $\Omega\in H^2(X)((\hbar))$ is the characteristic class of the star product $\star$ on $\calA^\hbar$.
\end{theorem}

As the pairing between the cyclic homology and cohomology of
$\calA^{\hbar}\rtimes \sfG$ is nondegenerate, we are able to state a
corollary about the Chern character of the element
$\widehat{P(E)}-\widehat{P(F)}$ in $K_0(\calA^{\hbar}\rtimes \sfG)$.
\begin{corollary}\label{cor:chern}
Let $\sfG$ be a proper \'etale groupoid. The image of the Chern character $Ch_k(\widehat{P(E)}-\widehat{P(F}))$ in $H^\bullet(\widetilde{X})((\hbar))$ under the map $\Psi$ (as is defined in Thm. \ref{chain-mor}) is equal to
\[
\frac{\hat{A}(B_0/\sfG) \, \Ch_\theta(\iota^*E-\iota^*F) \,
  \exp(-\frac{\iota^*\Omega}{2\pi \sqrt{-1}\hbar})}{\Ch_\theta(\lambda_{-1}N)},
\]
with the same notations as Theorem \ref{thm:orb-alg-index}.
\end{corollary}
\subsection{An index theorem for a constant Dirac structure on a 3-dim torus}
\label{Sec:IndConstDirac} In this subsection, we look at an
interesting example which is not a proper groupoid.

We consider a 3-dim torus $T^3=S^1\times S^1\times S^1$, and a
constant Dirac structure on $T^3$ defined by the following linear
Dirac structure on $\mathbb {R}^3$,
\[
{\mathcal {D}}_{\alpha, \beta, \theta}=\text{span} \left\{(\alpha, \beta,1;0,0,0),
(0,-\theta,0; 1,0,-\alpha), (\theta,0, 0; 0,1,-\beta)\right\}\subset
\mathbb {R}^3\times {\mathbb{R}^{3}}^*.
\]
We remark that $\alpha, \beta, \theta$ in the above definition of ${\mathcal {D}}_{\alpha, \beta, \theta}$ are all constant, and that ${\mathcal {D}}_{\alpha, \beta, \theta}$ defines a constant Dirac structure on $T^3$ with the identification $TT^3\oplus T^*T^3\simeq T^3\times \mathbb {R}^3\times {\mathbb{R}^{3}}^*$.
We fix $(\theta^1, \theta^2, \theta^3)$ to be the coordinates on
$T^3$. The Dirac structure ${\mathcal {D}}_{\alpha, \beta, \theta}$ defines a
characteristic foliation, which is spanned by the vector field
$Z=\alpha\partial_{\theta^1}+\beta\partial_{\theta^2}+\partial_{\theta^3}$.
We assume that both $\alpha$ and $\beta$ are irrational numbers.
Accordingly, the flow generated by the vector field $Z$ is ergodic
on $T^3$.

We consider the foliation ${\mathcal {F}}$ generated by $Z$ and
choose a complete transversal $X=\{\theta^3=0\}$ to ${\mathcal
{F}}$, which is a 2 dimensional torus. The holonomy group ${\mathbb
Z}$ acts on $X$ by translation mapping $(\theta^1, \theta^2)$ to
$(\theta^1+\alpha, \theta^2+\beta)$. We observe that the $\mathbb
{Z}$ action preserves the constant metric on $X$. Hence, the
holonomy groupoid $X\rtimes \mathbb {Z}$ satisfies the conditions
assumed in Section \ref{Sec:QuantConv}. Furthermore, as is
explained in \cite[Sec. 3]{tw:dirac}, the Dirac structure ${\mathcal
{D}}_{\alpha, \beta, \theta}$ defines a constant Poisson structure $\pi$ on $X$, i.e.
$\pi=\theta\partial_{\theta^1}\wedge
\partial_{\theta^2}$, which is invariant under the holonomy group action.
Therefore, we can apply the methods in Section \ref{Sec:QuantConv}
to construct a deformation quantization $\calA^\hbar_X\rtimes
\mathbb {Z}$ of the groupoid algebra $\calC^\infty(X)\rtimes
\mathbb{Z}$ (the transversal $X$ is a 2-dimensional torus, and there is no distinction between $\calC^\infty(X)$ and $\calC^\infty_c(X)$). In \cite[Sec. 3]{tw:dirac} the third author with
Weinstein used this idea to quantize an arbitrary constant Dirac
structure on an $n$-dimensional torus. Here, our construction is
slightly different as we are considering the formal deformation
quantization of the smooth algebra $\calC^\infty(X)\rtimes \mathbb {Z}$
instead of the $C^*$-algebras $C(X)\rtimes \mathbb {Z}$. We notice
that $X$ has a flat connection $d$ which is symplectic with respect
to the Poisson structure $\pi$ and also invariant under the $\Z$
action. Therefore, we choose the Fedosov connection in this example
to be
\[
D=d\theta^i\wedge (\partial_{\theta^i}-\partial_{y^i}),
\]
where $(y^1,y^2)$ are coordinates on the Weyl algebra bundle
$\calW$. In this case, the connection 1-form $A$ is equal to
$-1/\theta(y^1d\theta^2-y^2d\theta^1)$. The flat sections with
respect to the connection $D$ defines a deformation quantization
$\calA^\hbar_X$ of $\calC^\infty(X)$, and the crossed product algebra
$\calA^\hbar_X\rtimes \Z$ is a deformation quantization of
$\calC^\infty(X)\rtimes \Z$.

We apply Theorem \ref{chain-mor} to compute the cyclic homology of
the algebra $\calA^\hbar_X\rtimes \mathbb{Z}$. Notice that $\Z$
action on $X$ is free. Therefore, the Burghelea space $B_0$ of the groupoid $X\rtimes
\Z$ is equal to $X$, and the cyclic homology of
$\calA^\hbar_X\rtimes \Z$ is equal to the cyclic homology of the
groupoid $X\rtimes \Z$. As the spectral sequence associated to the
double complex $\Gamma_\cpt (X\rtimes \Z_\bullet,
\Omega^\bullet)$ degenerates at $E_2$, the homology of the groupoid
$X\rtimes Z$ is computed as follows
\[
\begin{split}
H_{1}(X\rtimes Z; \C((\hbar)))&=H_{-2}(X\rtimes Z;
\C((\hbar)))=\C((\hbar)),\\
H_0(X\rtimes Z; \C((\hbar)))&=H_{-1}(X\rtimes Z;
\C((\hbar)))=\C((\hbar))^{\times 3}.
\end{split}
\]
Similarly the cohomology of the groupoid $X\rtimes Z$ is computed by
\[
\begin{split}
H^0(X\rtimes Z; \C((\hbar)))&=H^3(X\rtimes Z;
\C((\hbar)))=\C((\hbar)),\\
H^1(X\rtimes Z; \C((\hbar)))&=H^2(X\rtimes Z;
\C((\hbar)))=\C((\hbar))^{\times 3}.
\end{split}
\]
We conclude with the following results about cyclic (co)homology of
$\calA^\hbar_X\rtimes \Z$.
\begin{corollary}
\label{cor:hom-tori}The cyclic homology and cohomology of
$\calA^\hbar_X\rtimes \Z$ is computed as follows,
\[
\begin{split}
HP_0(\calA^\hbar_X\rtimes \Z)=\C((\hbar))^{\times 4},\qquad
&HP_1(\calA^\hbar_X\rtimes \Z)=\C((\hbar))^{\times 4},\\
HP^0(\calA^\hbar_X\rtimes \Z)=\C((\hbar))^{\times 4},\qquad
&HP^1(\calA^\hbar_X\rtimes \Z)=\C((\hbar))^{\times 4}.
\end{split}
\]
\end{corollary}

We write out explicit cocycles generating the cohomology of the
cyclic vector space $\Tot^\bullet(\Omega^\bullet(B_\bullet))$
for the groupoid $X\rtimes \Z$. As $\Z$ acts on $X$ freely, $B_n$ is
a subspace of $X\times \Z^{\times (n+1)}$ isomorphic to
\[
\{(x, z_0, \cdots, z_n)|\quad x\in X, z_i\in \Z,  z_0+\cdots
+z_n=0.\}
\]
Let  $W^n$ be the delta function on $\Z$ supported
at $n$.
\begin{enumerate}
\item The degree 0 cohomology is generated by $\xi_0=1$ the constant
function.
\item The degree 1 cohomology is generated by $\eta_1=d\theta^1\in \Omega^1(X)^\Z$,
$\eta_2=d\theta^2\in \Omega^1(X)^\Z$ (we view $d\theta^1$ and
$d\theta^2$ as differential forms supported on the unit space of the
groupoid $X\rtimes \Z$), and $\eta_3=\delta=\sum_n nW^{-n}\otimes
W^{n}\in \Omega^0(B_1)$ (we view $ W^{-n}\otimes W^n$ a function on
$X\times \Z\times \Z$ taking 1 on $(x,-n,n)$ and 0 otherwise). We
remark that although $\sum_n nW^{-n}\otimes W^n$ is an infinite sum,
as a linear functional on $(\calC^\infty(X)\rtimes \Z)^{\otimes 2}$ it
is well defined because an element in $\calC^\infty(X)\rtimes \Z$ (and
also $\calA^\hbar_X\rtimes \Z$) has only finitely many components in
$\Z$.
\item The degree 2 cohomology is generated by $\xi_1=d\theta^1\wedge d\theta^2\in \Omega^2(X)^\Z$,
$\xi_2=d\theta^1\otimes \delta\in \Omega^1(B_1)^\Z$, and
$\xi_3=d\theta^2\otimes \delta\in \Omega^1(B_1)^\Z$, where $\delta$
is defined as above.
\item The degree 3 cohomology is generated by $\eta_0=d\theta^1\wedge d\theta^2 \otimes \delta\in \Omega^2(B_1)^\Z$.
\end{enumerate}

Next we compute explicit cyclic cocycles on $\calA^\hbar\rtimes \Z$
generating $HP^0(\calA_X^\hbar\rtimes \Z)$. (The same method can be
extended easily to compute cocycles in $HP^1(\calA^\hbar _X\rtimes
\Z)$.) As is explained before, the connection $A$ chosen for the
Fedosov connection defining $\calA^\hbar_X$ is equal to
$-1/\theta(y^1d\theta^2-y^2d\theta^1)$. This simplifies the
computation of the map $\Phi$ in Section 3.4. We find the following
generators in $HP^0(\calA_X^\hbar\rtimes \Z)$ by inserting the
connection $A$ into the definition of $\Phi$.
\begin{eqnarray*}
&\Phi(\xi_0)&(fW^n)=\\
&&\left\{\begin{array}{ll}-\frac{1}{\hbar\theta}\int_X f d\theta_1\wedge d\theta_2,&\qquad\qquad\qquad\qquad\qquad\qquad n=0\\
0&\qquad\qquad\qquad\qquad\qquad\qquad\text{otherwise}\end{array}\right.;\\
&\Phi(\xi_1)&(f_0W^{n_0}, f_1W^{n_1},f_2W^{n_2})=\\
&&\left\{\begin{array}{ll}\int_X \tau_2(\hat{f}_0,
\widehat{n_0^*(f_1)}, \widehat{n_0^*n_1^*(f_2)})d\theta^1\wedge
d\theta^2,&\qquad\qquad n_0+n_1+n_2=0\\ 0&\qquad\qquad\text{otherwise}\end{array}\right.;\\
&\Phi(\xi_2)&(f_0W^{n_0}, f_1W^{n_1}, f_2W^{n_2})=\\
&&\left\{\begin{array}{ll}\frac{\hbar}{\theta}\int_X n_1f_0\star
n_0^*(f_1)\star (-n_2)^*(\partial_{\theta^2}f_2)d\theta^1\wedge
d\theta^2,&\quad n_0+n_1+n_2=0,
\\ 0&\quad \text{otherwise}\end{array}\right.;\\
&\Phi(\xi_3)&(f_0W^{n_0}, f_1W^{n_1}, f_2W^{n_2})=\\
&&\left\{\begin{array}{ll}-\frac{\hbar}{\theta}\int_X n_1f_0\star
n_0^*(f_1)\star (-n_2)^*(\partial_{\theta^1}f_2)d\theta^1\wedge
d\theta^2,&n_0+n_1+n_2=0,
\\ 0&\text{otherwise.}\end{array}\right.
\end{eqnarray*}

We consider the following three subalgebras of $\calA^\hbar_X\rtimes \Z$.
\begin{enumerate}
\item The algebra $\calA^\hbar_X$ is embedded in $\calA^\hbar_X\rtimes
\Z$ as functions supported at the component of zero in $X\rtimes\Z$.
\item On $X$, we consider the subspace of functions constant along $\theta_2$. The
subspace of such functions is isomorphic to $\calC^\infty(S^1)$, the
algebra of smooth functions on $S^1$. The group $\Z$ action
preserves this subspace and therefore the associated crossed product
algebra defines a subalgebra of $\calA^\hbar_X\rtimes \Z$. This
subalgebra is isomorphic to the groupoid algebra
$\calC^\infty(S^1)\rtimes_\alpha \Z$ of the transformation groupoid that
$\Z$ acts on $S^1$ by translating angle $\alpha$.
\item Similarly to the previous case, we consider smooth functions
on $X$ constant along $\theta_1$, and obtain a subalgebra which
isomorphic to the groupoid algebra $\calC^\infty(S^1)\rtimes_\beta \Z$
of transformation groupoid of $\Z$ acting on $S^1$ by translating
angle $\beta$.
\end{enumerate}

We are interested in elements of
$K_0(\calA^\hbar_X\rtimes \Z)$ induced from $K_0$ elements associated to three subalgebras, $\calA^\hbar_X$, $\calC^\infty(S^1)\rtimes_\alpha \Z$, and
$\calC^\infty(S^1)\rtimes_\beta \Z$.

According to \cite[Sec. 6.1]{FedDQIT}, the $K_0(\calA^\hbar_X)$ is
isomorphic to $K_0(\calC^\infty(X))$. For an element $e$ in
$K_0(\calA_X^\hbar)\subset K_0(\calA^\hbar_X\rtimes \Z)$, we can
compute $Ch(e)$ directly by the results in \cite[Sec. 5]{ppt:high},
\[
Ch(e)=\widehat{A}(X)Ch(e)\exp(-\frac{\Omega}{\hbar}),
\]
where $Ch(e)$ is the Classical Chern character on $K_0(\calC^\infty(X))$ and
$\Omega$ is the characteristic classes of the invariant star product
$\star$ on $\calA^\hbar_X$, which in our case is $-1/\theta
d\theta^1\wedge d\theta^2$. This is a differential form supported on
the unit space of $X\rtimes \Z$.

As the consideration of subalgebras $\calC^\infty(S^1)\rtimes_\alpha\Z$
and $\calC^\infty(S^1)\rtimes_\beta \Z$ is a copy from one to the other,
we will focus on one of them, $\calC^\infty(S^1)\rtimes_\alpha \Z$. The
construction in \cite{rieffel:tori} defines two linear independent
elements in $K_0(\calC^\infty(S^1)\rtimes_\alpha \Z)\subset
K_0(\calA^\hbar_X\rtimes \Z)$,
\[
e_1=1,\qquad e_2=W^{-1}g+f+gW,
\]
where $f$ and $g$ are smooth functions defined on $S^1$ constructed
in \cite{rieffel:tori}. We recall the definition of $f$ and $g$.
Assume that $0\leq\alpha<1/2$ (the similar construction also works
for $\alpha\geq 1/2$), and choose $\epsilon>0$ such that
$\epsilon<\alpha$ and $\alpha+\epsilon<1/2$. The function $f$ is
defined to be any smooth monotone function on $[0,\epsilon]$ with
$f(0)=0$ and $f(\epsilon)=1$; on $[\alpha, \alpha+\epsilon]$, $f$ is
defined by $f(x)=1-f(x-\alpha)$; on $[\epsilon, \alpha]$ define $f$
to be 1; and on $[\alpha+\epsilon, 1]$ define $f$ to be 0. The
function $g$ is defined to be $(f-f^2)^{1/2}$ on $[\alpha,
\alpha+\epsilon]$ and zero everywhere else.

For the crossed product algebra ${\mathcal S}(\Z,
\calC^\infty(S^1))$ of Schwartz functions on $\Z$ valued in
$\calC^\infty(S^1)$, Pimsner-Voiculescu's exact sequence proves that
$K_0({\mathcal S}(\Z, \calC^\infty(S^1)))$ is isomorphic to
$\Z\oplus\Z$ with generators $e_1$ and $e_2$. And we can easily
check that all our cyclic cocycles  extend to be well defined on
this enlarged algebra. However, for our algebra which consists of
finitely supported functions on $\Z$ with value $\calC^\infty(S^1)$,
we do not have a full description of
$K_0(\calC^\infty(S^1)\rtimes_\alpha \Z)$, though from the
computation below we know that $e_1$ and $e_2$ are different
elements in $K_0(\calC^\infty(S^1)\rtimes _\alpha \Z)$.

We consider the evaluations of $\Phi(\xi_i)$, $i=1,\cdots, 4$ on
$e_1$ and $e_2$. As $e_1$ and $e_2$ are constant along $\theta_2$,
$\Phi(\xi_1)(e_1)=\Phi(\xi_1)(e_2)=\Phi(\xi_2)(e_1)=\Phi(\xi_2)(e_2)=0$.
The cocycle $\Phi(\xi_1)$ is $-1/\hbar\theta$ times the standard
trace on $\calC^\infty(S^1)\rtimes_\alpha \Z$. Therefore, we have
\[
\begin{split}
\Phi(\xi_0)(e_1)=-\frac{1}{\hbar\theta},&\qquad
\Phi(\xi_0)(e_2)=-\frac{\alpha}{\hbar\theta}, \\
\Phi(\xi_3)(e_1)=0,&\qquad \Phi(\xi_3)(e_2)=6
\frac{\hbar}{\theta}\int_{S^1}g^2df(\theta^1)=\frac{\hbar}{\theta}.
\end{split}
\]
Therefore, we conclude
\[
\Psi(Ch(e_1))=-\frac{1}{\hbar\theta},\qquad
\Psi(Ch(e_2))=-\frac{\alpha}{\hbar\theta}|_{id}+\frac{\hbar}{2\theta}d\theta_1|_{(1,-1)}
-\frac{\hbar}{2\theta}d\theta_1|_{(-1,1)},
\]
where $(\cdots)|_{(a,b)}$ means a differential form supported at the
$(a,b)$ component of $B(X\rtimes \Z)_2$.

\begin{remark}
We end this section with two remarks.
\begin{enumerate}
\item Our cocycle $\Phi(\xi_3)$ is closely related to the Chern character $c_1$ on
the quantum torus algebra ${\mathcal {S}}(\Z, \calC^\infty(S^1))$
introduced by Connes \cite{connes:foliation}; it differs by a scaler due to a
different normalization.
\item This subsection's discussion about algebraic index theory of a constant Dirac structure on a torus
can be generalized to arbitrary constant Dirac structures on an
arbitrary dimensional torus following the ideas developed in
\cite{tw:dirac} and the constructions in Section \ref{Sec:ExpChaMap}
of this paper. For a general constant Dirac structure on $T^n$, we
can choose a complete transversal $X$ to the characteristic
foliation as in \cite{tw:dirac} and consider the monodromy groupoid
$\sfG_X$. The unit space $\sfG_0$ of the groupoid $\sfG$ has a
$\sfG$ invariant metric, and the Dirac structure defines on $\sfG_0$
a $\sfG$ invariant regular Poisson structure. All considerations in
this paper can be directly generalized to such $\sfG_X$, and we will
leave the details to the readers.
\end{enumerate}
\end{remark}

\appendix
\bibliographystyle{alpha}

\begin{thebibliography}{}

\bibitem[{\sc BrGoNeTs}]{bgnt}
{\sc P.~Bressler, A.~Gorokhovsky, R.~Nest}, and {\sc B.~Tsygan},
\textit{Deformations of algebroid stacks}. \texttt{arXiv:0810.0030}
(2008).

\bibitem[{\sc Br}]{bry}
  {\sc J.L.~Brylinski}: {\it A differential complex for Poisson manifolds},
  J. Differential Geom.~\textbf{28}, no. 1, 93--114 (1988).

\bibitem[{\sc BrNi}]{bn} {\sc J.L.~Brylinski} and {\sc V.~Nistor}:
  \textit{Cyclic cohomology of \'{e}tale groupoids}.
  $K$-theory \textbf{8}, 341--365 (1994).

\bibitem[{\sc Ca}]{Sil:LSG}
   {\sc A.~Cannas da Silva}:
   \textit{Lectures on Symplectic Geometry}. Lec.~Notes Math.~\textbf{1764}, Springer (2001).


\bibitem[{\sc Co}94]{c:book} {\sc A.~Connes}:
  \textit{Noncommutative Geometry}. Academic Press, 1994.

\bibitem[{\sc Co}82]{connes:foliation} {\sc A.~Connes}:
  \textit{A survey of foliations and operator algebras}.
  Operator algebras and applications, Part I (Kingston, Ont., 1980), pp. 521--628,
  Proc. Sympos. Pure Math., \textbf{38}, Amer. Math. Soc., Providence, R.I., 1982.



\bibitem[{\sc Cr}]{crainic} {\sc M.~Crainic}:
  \textit{Cyclic cohomology of \'{e}tale groupoids: the general case}.
  $K$-theory \textbf{17}, 319--362. (1999).


\bibitem[{\sc Fe95}]{FedDQIT} {\sc B. Fedosov}:
  \textit{Deformation quantization and index theory}. Akademie Verlag, 1995.

\bibitem[{\sc Fe02}]{fe:g-index} {\sc B.~Fedosov}:
  \textit{On $G$-trace and $G$-index in deformation quantization}.
  Conference Mosh\'e Flato 1999 (Dijon), Lett. Math. Phys.~\textbf{52} ,
  no. 1, 29--49 (2002).

\bibitem[{\sc FeFeSh}]{ffs}
  {\sc B.~Feigin, G.~Felder}, and {\sc B.~Shoikhet}:
  \textit{Hochschild cohomology of the Weyl algebra and traces in
  deformation quantization}. Duke Math. J.~\textbf{127}, no. 3, 487--517 (2005).

\bibitem[{\sc GeJo}]{gj} {\sc E.~Getzler}, and {\sc J.D.~Jones}:
  \textit{The cyclic homology of crossed product algebras.}
  J. Reine Angew. Math.~\textbf{445} 161--174 (1993).


\bibitem[{\sc KhRa}]{kr} {\sc M.~Khalkhali}, and {\sc B.~Rangipour}:
  \textit{On the generalized cyclic Eilenberg-Zilber theorem.}
  Canad. Math. Bull.~ \textbf{47}, no. 1, 38--48 (2004).

\bibitem[{\sc Lo}]{loday}
  {\sc J.~Loday}:
  \textit{Cyclic Homology}, Grundlehren der mathematischen
  Wissenschaften Vol.~\textbf{301}, Springer, 1998.




\bibitem[{\sc Mo}]{molino:book}
  {\sc P.~Molino}: \textit{Riemannian Foliations}.
  Birkh\"auser, Boston (1988).

\bibitem[{\sc NeTs}95]{nets95} {\sc R.~Nest} and {\sc B.~Tsygan}:
  {\it Algebraic index theorem},
  Comm. Math. Phys {\bf 172}, 223--262 (1995).


\bibitem[{\sc NePfPoTa}]{nppt}
  {\sc N.~Neumaier, M.J.~Pflaum, H.~Posthuma}, and {\sc X.~Tang}:
  \textit{Homology of formal deformations of proper \'etale Lie
  groupoids}. J. Reine Angew. Math.~\textbf{593}, 117-168 (2006).



\bibitem[{\sc PfPoTa}07]{ppt:index}
  {\sc M.J.~Pflaum, H.~Posthuma}, and {\sc X.~Tang}:
  \textit{An algebraic index theorem for orbifolds}.
  Adv. Math.~\textbf{210}, p. 83--121 (2007).

\bibitem[{\sc PfPoTa}08]{ppt:high}
  {\sc M.J.~Pflaum, H.~Posthuma}, and {\sc X.~Tang}:
  \textit{Cyclic cocycles on deformation quantizations and higher index theorems}.
  \texttt{arXiv:0805.1411} (2008).

\bibitem[{\sc PfPoTaTs}]{pptt}
  {\sc M.J.~Pflaum, H.~Posthuma, X.~Tang} and {\sc H.-H.~Tseng}:
  \textit{Orbifold cup products and ring structures on Hochschild
  cohomologies}.
  \texttt{arXiv:0706.0027} (2007).

\bibitem[{\sc Ph}]{philips:k-theory}
  {\sc N.Ch.~Phillips}:
  \textit{Equivariant $K$-theory for proper actions}.
  Pitman Research Notes in Mathematics Series, \textbf{178}. Longman Scientific {\&}
  Technical, Harlow; copublished in the United States with John Wiley {\&}   Sons, Inc.,
  New York, 1989.

\bibitem[{\sc Ri}]{rieffel:tori}
  {\sc M.A.~Rieffel}:
  \textit{$C\sp{\ast} $-algebras associated with irrational rotations}.
  Pacific J. Math. \textbf{93}, no. 2, 415--429 (1981.


\bibitem[{\sc Ta}06]{ta}
  {\sc X.~Tang}:
  \textit{Deformation quantization of pseudo Poisson
  groupoids}, Geom.~Func.~Analysis \textbf{16} Nr.~3, 731--766 (2006).

\bibitem[{\sc TaWe}]{tw:dirac}
  {\sc X.~Tang} and {\sc A. Weinstein}:
  \textit{Quantization and Morita equivalence for constant Dirac structures on tori}.
  Ann. Inst. Fourier (Grenoble) \textbf{54}, no. 5, 1565--1580, xvi, xxii (2004).

\bibitem[{\sc We}]{weibel}
  {\sc Ch.~Weibel}:
  \textit{An Introduction to Homological Algebra},
  Cambridge Studies in Advanced Mathematics \textbf{38},
  Cambridge Univ. Press 1995.
\end{thebibliography}

\end{document}